\documentclass[final, 12pt,times]{elsarticle}

\usepackage{url}
\usepackage[utf8]{inputenc}
\usepackage[T1]{fontenc}
\usepackage{makecell} 
\usepackage{booktabs} 
\usepackage{array}    
\usepackage{placeins}

\usepackage{subcaption}
\usepackage{algpseudocode}
\usepackage{algorithm}
\usepackage{fancyhdr}
\usepackage{comment}
\usepackage{xcolor}
\usepackage{bm}
\usepackage[pagewise]{lineno}

\usepackage{tikz,graphicx}
\usepackage{multirow}
\usepackage{amsfonts,amsmath,amssymb}
\usepackage{mathtools}
\usepackage{mathrsfs}
\usepackage{upgreek}

\begin{document}

\begin{frontmatter}

\title{Data-Assimilation-Assisted Reinforcement Learning\\
for Power Grid Control under Load Uncertainty}

\author[au]{Ruoyu Hu\corref{cor1}}
\ead{ruh0003@auburn.edu}

\author[math]{Feng Bao}
\ead{fbao@fsu.edu}

\author[au]{Zezhong Zhang}
\ead{zez0002@auburn.edu}

\author[au]{Yanzhao Cao}
\ead{yzc0009@auburn.edu}

\author[ornl]{Guannan Zhang}
\ead{zhangg@ornl.gov}

\cortext[cor1]{Corresponding author}

\affiliation[au]{
    organization={Department of Mathematics and Statistics, Auburn University},
    city={Auburn},
    postcode={36849},
    state={Alabama},
    country={USA}
}

\affiliation[math]{
    organization={Department of Mathematics, Florida State University},
    city={Tallahassee},
    postcode={32306},
    state={Florida},
    country={USA}
}

\affiliation[ornl]{
    organization={Computer Science and Mathematics Division, Oak Ridge National Laboratory},
    city={Oak Ridge},
    postcode={37831},
    state={Tennessee},
    country={USA}
}

\begin{abstract}
Reliable power grid control requires sequential decisions under transmission constraints,
time-varying demand, renewable variability, and imperfect load information. We investigate
how load information quality affects control in a customized transmission network simulator
inspired by Grid2Op and based on a DC power flow model. The environment includes stochastic
spatial loads, temporally correlated renewable-like generation, energy storage, line
protection, generator dispatch, and line reconnection. An Ensemble Score Filter (EnSF)
corrects noisy load information before it is passed to the controller. A heuristic study
first isolates the effect of this correction under a fixed feedback rule. Proximal policy
optimization (PPO) agents are then trained for sequential generator and reconnection control,
frozen, and evaluated with Forward, EnSF-corrected, and Truth load information on paired
stochastic scenarios. EnSF consistently lowers load estimation error. Under heuristic
control, the correction extends survival and reduces warning and overload exposure. With
frozen PPO policies, it increases cumulative return and keeps performance closer to the
Truth benchmark, with a larger advantage when load uncertainty is amplified beyond the
nominal training distribution. In both experiments, EnSF reduces load-estimation error and improves the resulting control performance in the simplified transmission network.
\end{abstract}

\end{frontmatter}

\section{Introduction}
\label{sec:introduction}

Electric power systems must continuously balance generation and demand while respecting
transmission and equipment limits. The problem becomes more difficult as demand varies
over time and larger amounts of renewable generation and storage enter the system.
Changes in injections at one location can alter flows throughout the network, and actions
that relieve current stress may reduce flexibility later. Renewable variability, uncertain
demand, line contingencies, and an expanding set of control actions make power grid
operation a challenging sequential decision problem
\cite{glavic2019,marot2021,vandersar2025}.

Control decisions also depend on the quality of the information available at the time of
action. Load forecasts, renewable forecasts, and network measurements generally differ
from the physical quantities they represent. An error in estimated demand can alter
generation response, storage charging or discharging, and transmission flows. Network
coupling allows a local information error to affect operating constraints elsewhere, so
robustness to perturbed observations is an important issue for learning-based power system
control \cite{pan2021}.

Reinforcement learning has been increasingly used to study sequential control problems
in power systems \cite{glavic2019}. Grid2Op provides a simulation platform for repeated
decision making in transmission networks, while the Learning to Run a Power Network
(L2RPN) challenges established common benchmarks for learning-based grid control
\cite{developers2026,marot2021}. A central difficulty in these problems is the large and
highly structured action space associated with generator operation, line switching, and
network topology. Hierarchical policies and afterstate representations have been
introduced to reduce this complexity, while other approaches use restricted action sets,
simulation-based action screening, or combinations of reinforcement learning and
rule-based control \cite{yoon2021,zhou2021,chauhan2023,lehna2023}. Graph-based and
distributed policies have also been investigated to improve scalability
\cite{fabrizio2025}. Most of these studies concentrate on the control side of the problem,
including how actions are represented, selected, and constrained. Uncertainty in the
state or injection information supplied to the controller is often absorbed into the
environment rather than examined separately.

Data assimilation provides a way to treat this information uncertainty explicitly. A
forward model predicts the evolving state, while observations provide additional
information about the physical system. Sequential data assimilation combines these
sources to approximate the conditional state distribution and update the state estimate
as new observations become available \cite{carrassi2018}. The Kalman filter gives the
optimal recursive estimator for linear Gaussian systems \cite{kalman1960,kalman1961}.
Ensemble Kalman filters (EnKF) extend this framework to large nonlinear systems through
ensemble approximations and have become widely used because of their computational
efficiency \cite{evensen1994,houtekamer1998,houtekamer2016}. Their updates, however,
are based on Gaussian-type representations through ensemble means and covariances,
which can become restrictive when the filtering distribution or the observation operator
is strongly nonlinear or non-Gaussian.

Particle filters provide a more general Bayesian representation through weighted
particles \cite{gordon1993,vanleeuwen2019}. They can represent non-Gaussian filtering
distributions, but the particle weights tend to concentrate as the dimension increases,
leading to weight degeneracy and rapidly growing ensemble requirements. Nonlinear
filtering can also be formulated through stochastic partial differential equations for
conditional densities, including the Zakai equation and related numerical methods
\cite{gobet2006,hu2002,bao2014,bao2018,bao2020FbdsdeFilter}. These formulations are
general, but their numerical cost can become substantial in high-dimensional systems.
High dimensionality and nonlinear observations remain major difficulties in
nonlinear filtering.

Score-based filtering offers an alternative representation of high-dimensional,
non-Gaussian filtering distributions. The score-based filter characterizes the filtering
density through its logarithmic density gradient and uses a diffusion process to generate
samples from that density \cite{bao2024ScoreBased}. The Ensemble Score Filter (EnSF)
avoids training a neural score model at each filtering stage. Instead, the score is
approximated directly from forecast ensembles using a closed form representation and a
Monte Carlo approximation, while observational information enters through a likelihood
score update \cite{bao2024EnsembleScore}. Avoiding repeated score network training is
particularly useful in sequential data assimilation, where filtering must be performed
at every observation time. EnSF has been demonstrated for high-dimensional nonlinear
systems and has subsequently been applied to porous media flow, partial observation
problems, observation network studies, earth system prediction, and real energy consumption forecasting
\cite{hu2026TwoPhaseFlow,liang2025,xiong2026,hu2026EnergyForecast,iensf, wang2026global}.

An imperfect load forecast can directly affect the next grid control action. Without
filtering, new measurements do not systematically correct that forecast before generator
and storage decisions are made, so estimation error can propagate into dispatch, storage
use, and transmission flows. Data assimilation inserts a correction between prediction
and control without changing the downstream control law. Filtering has improved feedback
and reinforcement learning decisions under noisy measurements \cite{song2026}, but
high-dimensional nonlinear load correction has received much less direct attention in
sequential transmission grid control.

We treat load information as the uncertain state component supplied to the grid
controller and use EnSF for nonlinear data assimilation. The experiments separate load
correction from controller design. The first applies one heuristic feedback rule under
Forward, EnSF-corrected, and Truth load information, so only the information entering the
controller changes. The second trains proximal policy optimization (PPO) policies for
sequential generator and line reconnection control and freezes them before changing the
load information channel. This allows the same learned policy to be evaluated with
different load estimates. Paired physical scenarios are used throughout so that each
information mode sees the same stochastic grid realization.

The remainder of the paper is organized as follows. Section~\ref{sec:grid_control_problem}
introduces the power grid control problem and its relation to Grid2Op and L2RPN.
Section~\ref{sec:data_assimilation} presents the Bayesian filtering formulation and the
Ensemble Score Filter. Section~\ref{sec:numerical_implementation} describes the
numerical power grid and filtering implementation. Section~\ref{sec:numerical_examples}
presents the heuristic and reinforcement learning experiments. The final section
summarizes the results, discusses the limitations of the simplified grid model, and
outlines directions for future work.

\section{Power Grid Control Problem}
\label{sec:grid_control_problem}

We consider sequential control of a transmission network represented by a graph
\[
    \mathcal G=(\mathcal V,\mathcal E),
\]
where $\mathcal V$ is the set of substations and $\mathcal E$ is the set of transmission
lines. At time $t_n$, the active power injection at substations $i\in\mathcal V$ is determined
by the generation, renewable production, storage power, and load connected to that substation.
Let $\mathcal G_i$, $\mathcal R_i$, and $\mathcal D_i$ denote the sets of controllable
generators, renewable generators, and loads connected to substation $i$, respectively. Using
the sign convention that positive storage power corresponds to discharge into the grid,
the net injection is
\begin{equation}
    p_{i,n}
    =
    \sum_{g\in\mathcal G_i} P^g_{g,n}
    +
    \sum_{r\in\mathcal R_i} P^r_{r,n}
    +
    P^s_{i,n}
    -
    \sum_{d\in\mathcal D_i} P^d_{d,n},
    \label{eq:bus_injection}
\end{equation}
where $P^g_{g,n}$ is the active power output of controllable generator $g$,
$P^r_{r,n}$ is the active power output of renewable generator $r$,
$P^s_{i,n}$ is the net storage injection at substation $i$, and $P^d_{d,n}$ is the
active power demand of load $d$. The active power injections satisfy the balance condition over the full system
\begin{equation}
    \sum_{i\in\mathcal V} p_{i,n}=0.
    \label{eq:power_balance}
\end{equation}
Generation and storage are limited by their operating ranges, while load and renewable
generation vary with time and act as exogenous inputs to the grid.

Active power transmission is modeled using a direct-current (DC) power flow
approximation. Let $\theta_{i,n}$ denote the voltage phase angle at substation $i$ at time
$t_n$. For a transmission line $\ell=(i,j)\in\mathcal E$, let $b_\ell$ denote its
line susceptance. The active power flow on the line is approximated by
\begin{equation}
    F_{\ell,n}
    =
    b_\ell\left(\theta_{i,n}-\theta_{j,n}\right),
    \label{eq:dc_line_flow}
\end{equation}
where $F_{\ell,n}$ is the active power flow from substation $i$ to substation $j$ under the selected
line orientation. Together with the nodal power balance equations, these relations
determine the line flows for a given set of injections and line statuses. The phase angle
of one reference substation is fixed to remove the arbitrary angular offset. The DC formulation
retains the dependence of transmission flows on spatially distributed active power
injections while neglecting voltage magnitude and reactive power effects.

For each transmission line $\ell$, let $\overline F_\ell$ denote its thermal
active power limit. The line loading ratio is defined as
\begin{equation}
    \rho_{\ell,n}
    =
    \frac{|F_{\ell,n}|}{\overline F_\ell}.
    \label{eq:line_loading_ratio}
\end{equation}
A line is overloaded when $\rho_{\ell,n}>1$, and the maximum line loading ratio at time
$t_n$ is
\begin{equation}
    \rho_n^{\max}
    =
    \max_{\ell\in\mathcal E}\rho_{\ell,n}.
    \label{eq:max_line_loading}
\end{equation}
Lines that remain overloaded for an extended period can trip because of thermal stress
or the activation of protection mechanisms. The loss of a transmission line changes the
network configuration and redistributes power over the remaining lines. This
redistribution can increase loading elsewhere in the network and, in severe cases, lead
to additional line trips or loss of feasible operation.

Energy storage provides a controllable means of absorbing or supplying active power.
Let $E_n$ denote the stored energy at time $t_n$, and let
$\Delta t=t_{n+1}-t_n$ denote the physical time step. We denote charging and discharging
powers by $P_n^{\rm ch}\geq0$ and $P_n^{\rm dis}\geq0$, respectively. With charging
efficiency $\eta_{\rm ch}\in(0,1]$ and discharging efficiency
$\eta_{\rm dis}\in(0,1]$, the storage energy evolves according to
\begin{equation}
    E_{n+1}
    =
    E_n
    +
    \Delta t
    \left(
        \eta_{\rm ch} P_n^{\rm ch}
        -
        \frac{P_n^{\rm dis}}{\eta_{\rm dis}}
    \right).
    \label{eq:storage_energy}
\end{equation}
The storage variables satisfy
\begin{equation}
    0\leq E_n\leq \overline E,
    \qquad
    0\leq P_n^{\rm ch}\leq \overline P^{\rm ch},
    \qquad
    0\leq P_n^{\rm dis}\leq \overline P^{\rm dis},
    \label{eq:storage_constraints}
\end{equation}
where $\overline E$ is the storage energy capacity and
$\overline P^{\rm ch}$ and $\overline P^{\rm dis}$ are the maximum charging and
discharging powers. The corresponding net storage injection is
\begin{equation}
    P_n^s=P_n^{\rm dis}-P_n^{\rm ch},
    \label{eq:storage_net_power}
\end{equation}
so that $P_n^s>0$ represents net discharge into the grid and $P_n^s<0$ represents
net charging. Storage couples decisions across time. Discharging can reduce
an instantaneous generation deficit but also reduces the energy available at later
times, while charging has the opposite effect.

Although each electrical operating point is obtained from an algebraic power flow
problem, grid operation remains sequential. Load and renewable generation change with
time, storage carries energy between time steps, and line trips and reconnection
restrictions alter the network available to later decisions. The controller acts on the
current network condition together with the operational state inherited from previous
steps. We use generator output adjustment and line reconnection as the principal control
variables.

Grid2Op was developed as a platform for studying this type of sequential decision making
in power systems \cite{donnot2020}. The Learning to Run a Power Network (L2RPN)
initiative extended this setting through a series of power grid control
competitions in which agents were evaluated under common network models, operating
scenarios, and evaluation rules \cite{marot2020,marot2021}. These competitions emphasized
the ability of an agent to maintain transmission network operation over long time
horizons under changing injections and network conditions. A major difficulty is the
large combinatorial action space, especially when line switching and busbar topology
changes are available. Competition entries and subsequent studies have used
hierarchical policies, action screening, search, and combinations of rule-based and
learning-based control to manage this complexity.

The model used here is a simplified version of the broader Grid2Op/L2RPN setting. It
retains spatially distributed generation and demand, transmission constraints, storage,
line outages, and reconnection, but omits the full busbar topology action space. The
smaller action space makes it easier to isolate uncertain state information and track how
data assimilation affects subsequent grid operation.

Let
\[
    P_n^d
    =
    \left(P^d_{1,n},\ldots,P^d_{N_d,n}\right)
\]
denote the true load vector at time $t_n$, where $N_d$ is the number of load
components in the network. This vector determines the physical load applied to the
power flow problem. The controller, however, may act on an estimated load vector
\[
    \widehat P_n^d
    =
    \left(\widehat P^d_{1,n},\ldots,\widehat P^d_{N_d,n}\right).
\]
A discrepancy between $P_n^d$ and $\widehat P_n^d$ can change generator response,
storage use, and transmission flows. Estimating $\widehat P_n^d$ is thus a filtering
problem in which observations are used to correct model predictions.

\section{Data Assimilation and the Ensemble Score Filter}
\label{sec:data_assimilation}

\subsection{Data assimilation problem and Bayesian filtering framework}
\label{sec:da_bayesian_filtering}

Let $X_{t_n}\in\mathbb R^{d_x}$ denote the state of a dynamical system at time $t_n$,
where $d_x$ is the state dimension. We consider the state model
\begin{equation}
    X_{t_{n+1}}
    =
    f_n\left(X_{t_n},\omega_{t_n}\right),
    \qquad
    n=0,1,2,\ldots,
    \label{eq:da_state}
\end{equation}
where $f_n$ is the forward model and $\omega_{t_n}$ represents model uncertainty or
other stochastic inputs. The observation model is
\begin{equation}
    Y_{t_{n+1}}
    =
    h_{n+1}\left(X_{t_{n+1}}\right)
    +
    \epsilon_{t_{n+1}},
    \label{eq:da_observation}
\end{equation}
where $Y_{t_{n+1}}\in\mathbb R^{d_y}$ is the observation vector, $d_y$ is the
observation dimension, $h_{n+1}:\mathbb R^{d_x}\rightarrow\mathbb R^{d_y}$ is the
observation operator, and $\epsilon_{t_{n+1}}$ denotes observational noise. We assume
Gaussian observational noise
\begin{equation}
    \epsilon_{t_{n+1}}
    \sim
    \mathcal N(0,R_{n+1}),
    \label{eq:da_observation_noise}
\end{equation}
where $R_{n+1}\in\mathbb R^{d_y\times d_y}$ is the observation error covariance
matrix.

Let
\begin{equation}
    \mathcal Y_n
    :=
    \sigma\left(
        Y_{t_1},\ldots,Y_{t_n}
    \right)
    \label{eq:da_observation_filtration}
\end{equation}
denote the $\sigma$-algebra generated by the observations available up to time $t_n$.
The optimal filtering estimate of the state at $t_{n+1}$ is the conditional expectation
\begin{equation}
    \widehat X_{t_{n+1}}
    =
    \mathbb E
    \left[
        X_{t_{n+1}}
        \mid
        \mathcal Y_{n+1}
    \right].
    \label{eq:da_optimal_estimate}
\end{equation}
For nonlinear systems, the conditional distribution is generally not available in closed
form. We use the filtering density
\begin{equation}
    \pi_{n|n}(x)
    :=
    p\left(
        X_{t_n}=x
        \mid
        \mathcal Y_n
    \right),
    \label{eq:da_filtering_density}
\end{equation}
which characterizes the uncertainty in the state after assimilating observations through
time $t_n$.

Bayesian filtering proceeds recursively through prediction and update. Given the posterior
density $\pi_{n|n}$, the Chapman--Kolmogorov equation gives the prior density at
$t_{n+1}$ as
\begin{equation}
    \pi_{n+1|n}(x)=\int_{\mathbb R^{d_x}}p\left(X_{t_{n+1}}=x\mid X_{t_n}=x'\right)\pi_{n|n}(x')\,dx'.
    \label{eq:da_prediction}
\end{equation}
Here, $p(X_{t_{n+1}}=x\mid X_{t_n}=x')$ is the transition density induced by the
state model in Eq.~\eqref{eq:da_state}. The density $\pi_{n+1|n}$ contains the
information propagated by the forward model together with all observations assimilated
through time $t_n$.

After $Y_{t_{n+1}}$ becomes available, define the likelihood function
\begin{equation}
    \ell_{n+1}(x)
    :=
    p\left(
        Y_{t_{n+1}}
        \mid
        X_{t_{n+1}}=x
    \right).
    \label{eq:da_likelihood_definition}
\end{equation}
Bayes' rule gives the posterior filtering density
\begin{equation}
    \pi_{n+1|n+1}(x)=\frac{\ell_{n+1}(x)\pi_{n+1|n}(x)}{\displaystyle\int_{\mathbb R^{d_x}}\ell_{n+1}(x')\pi_{n+1|n}(x')\,dx'}.
    \label{eq:da_update}
\end{equation}
For the Gaussian observation model, the likelihood has the form
\begin{equation}
    \ell_{n+1}(x)\propto\exp\!\left[-\frac{1}{2}\left(h_{n+1}(x)-Y_{t_{n+1}}\right)^\top R_{n+1}^{-1}\left(h_{n+1}(x)-Y_{t_{n+1}}\right)\right].
    \label{eq:da_likelihood}
\end{equation}

\subsection{Ensemble Score Filter}
\label{sec:ensf}

The Ensemble Score Filter (EnSF) provides an ensemble-based approximation of the
prior and posterior filtering distributions introduced above
\cite{bao2024ScoreBased,bao2024EnsembleScore}. Its construction uses a score-based
diffusion process to represent a probability distribution through its logarithmic
density gradient and to generate samples from that distribution without training a
neural score model.

Let $Z_\tau\in\mathbb R^{d_x}$ denote an auxiliary diffusion process indexed by the
pseudo-time variable $\tau$. The pseudo-time is independent of the physical time
$t_n$ in the filtering problem. We consider the forward stochastic differential
equation
\begin{equation}
    dZ_\tau
    =
    b(\tau)Z_\tau\,d\tau
    +
    \sigma(\tau)\,dW_\tau,
    \label{eq:ensf_forward_sde}
\end{equation}
where $W_\tau$ is a $d_x$-dimensional standard Brownian motion,
$b(\tau)$ is the drift coefficient, and $\sigma(\tau)$ is the diffusion coefficient.
The coefficients are defined by
\begin{equation}
    b(\tau)
    =
    \frac{d\log\alpha_\tau}{d\tau},
    \qquad
    \sigma^2(\tau)
    =
    \frac{d\beta_\tau^2}{d\tau}
    -
    2\frac{d\log\alpha_\tau}{d\tau}\beta_\tau^2,
    \label{eq:ensf_coefficients}
\end{equation}
with
\begin{equation}
    \alpha_\tau=1-\tau,
    \qquad
    \beta_\tau^2=\tau.
    \label{eq:ensf_schedule}
\end{equation}

Let $Q_0$ denote a target probability density at pseudo-time $\tau=0$. For a fixed
initial value $Z_0=z_0$, the transition density of the linear forward diffusion is
\begin{equation}
    Q_\tau(z\mid z_0)
    =
    \mathcal N
    \left(
        z;
        \alpha_\tau z_0,
        \beta_\tau^2 I_{d_x}
    \right),
    \label{eq:ensf_transition_density}
\end{equation}
where $I_{d_x}$ is the $d_x\times d_x$ identity matrix. As $\tau\rightarrow1$,
$\alpha_\tau\rightarrow0$ and $\beta_\tau^2\rightarrow1$, so the diffused distribution
approaches the standard Gaussian distribution.

Let $Q_\tau(z)$ denote the marginal density of $Z_\tau$. Its score function is defined by
\begin{equation}
    \mathcal S(z,\tau)
    :=
    \nabla_z\log Q_\tau(z).
    \label{eq:ensf_score}
\end{equation}
Given this score, samples from $Q_0$ can be generated by solving the reverse-time
stochastic differential equation
\begin{equation}
    dZ_\tau
    =
    \left[
        b(\tau)Z_\tau
        -
        \sigma^2(\tau)
        \mathcal S(Z_\tau,\tau)
    \right]d\tau
    +
    \sigma(\tau)d\overleftarrow W_\tau,
    \label{eq:ensf_reverse_sde}
\end{equation}
where $\overleftarrow W_\tau$ denotes the Brownian motion associated with the
reverse-time process. Starting from the Gaussian reference distribution, the
reverse-time diffusion recovers samples from the target distribution as
$\tau$ decreases toward zero.

The score in Eq.~\eqref{eq:ensf_score} can be expressed directly in terms of the
target density. Since
\begin{equation}
    Q_\tau(z)
    =
    \int_{\mathbb R^{d_x}}
    Q_\tau(z\mid z_0)Q_0(z_0)\,dz_0,
    \label{eq:ensf_marginal_density}
\end{equation}
differentiation with respect to $z$ gives
\begin{equation}
    \mathcal S(z,\tau)
    =
    \int_{\mathbb R^{d_x}}
    -
    \frac{z-\alpha_\tau z_0}{\beta_\tau^2}
    w_\tau(z,z_0)
    Q_0(z_0)\,dz_0,
    \label{eq:ensf_score_integral}
\end{equation}
where
\begin{equation}
    w_\tau(z,z_0)
    =
    \frac{
        Q_\tau(z\mid z_0)
    }{
        \displaystyle
        \int_{\mathbb R^{d_x}}
        Q_\tau(z\mid z'_0)
        Q_0(z'_0)\,dz'_0
    }.
    \label{eq:ensf_weight}
\end{equation}
For each fixed $(z,\tau)$, the weight satisfies
\begin{equation}
    \int_{\mathbb R^{d_x}}
    w_\tau(z,z_0)Q_0(z_0)\,dz_0
    =
    1.
    \label{eq:ensf_weight_normalization}
\end{equation}

EnSF approximates Eq.~\eqref{eq:ensf_score_integral} by Monte Carlo sampling.
Suppose that $\{z^{(m)}\}_{m=1}^{M}$ is an ensemble of $M$ samples from $Q_0$.
Using a mini-batch of $J\leq M$ ensemble members, the score is approximated by
\begin{equation}
    \overline{\mathcal S}(z,\tau)
    =
    \sum_{j=1}^{J}
    -
    \frac{z-\alpha_\tau z^{(j)}}{\beta_\tau^2}
    \overline w_\tau\left(z,z^{(j)}\right),
    \label{eq:ensf_mc_score}
\end{equation}
where
\begin{equation}
    \overline w_\tau\left(z,z^{(j)}\right)
    =
    \frac{
        Q_\tau\left(z\mid z^{(j)}\right)
    }{
        \displaystyle
        \sum_{\ell=1}^{J}
        Q_\tau\left(z\mid z^{(\ell)}\right)
    }.
    \label{eq:ensf_mc_weight}
\end{equation}
Here, $M$ is the full ensemble size and $J$ is the number of ensemble members used
in the Monte Carlo score approximation. In practice, $J$ may be much smaller than $M$.
The score is evaluated directly from the ensemble samples and the known transition
density $Q_\tau(z\mid z_0)$.

For sequential filtering, let
$\{x_{n|n}^{(m)}\}_{m=1}^{M}$ denote an ensemble representing the posterior
filtering density $\pi_{n|n}$. The prediction step propagates each ensemble member
through the state model,
\begin{equation}
    x_{n+1|n}^{(m)}
    =
    f_n
    \left(
        x_{n|n}^{(m)},
        \omega_{t_n}^{(m)}
    \right),
    \qquad
    m=1,\ldots,M,
    \label{eq:ensf_prediction}
\end{equation}
where $\omega_{t_n}^{(m)}$ is the realization of model uncertainty associated with
ensemble member $m$. The predicted ensemble
$\{x_{n+1|n}^{(m)}\}_{m=1}^{M}$ provides a sample approximation of the prior
filtering density $\pi_{n+1|n}$. Using this ensemble in
Eq.~\eqref{eq:ensf_mc_score} gives the prior score
\begin{equation}
    \overline{\mathcal S}_{n+1|n}(z,\tau).
    \label{eq:ensf_prior_score}
\end{equation}

The observation at $t_{n+1}$ is incorporated through the likelihood
$\ell_{n+1}$ defined in Eq.~\eqref{eq:da_likelihood_definition}. EnSF approximates
the posterior score by
\begin{equation}
    \overline{\mathcal S}_{n+1|n+1}(z,\tau)
    =
    \overline{\mathcal S}_{n+1|n}(z,\tau)
    +
    g(\tau)
    \nabla_z\log\ell_{n+1}(z),
    \label{eq:ensf_posterior_score}
\end{equation}
where $g(\tau)$ is a damping function controlling the contribution of the likelihood
along the pseudo-time interval. It satisfies
\begin{equation}
    g(0)=1,
    \qquad
    g(1)=0,
    \label{eq:ensf_damping_conditions}
\end{equation}
and we take
\begin{equation}
    g(\tau)=1-\tau.
    \label{eq:ensf_damping}
\end{equation}

The posterior score
$\overline{\mathcal S}_{n+1|n+1}$ is then used in the reverse-time diffusion
Eq.~\eqref{eq:ensf_reverse_sde} to generate an updated ensemble
$\{x_{n+1|n+1}^{(m)}\}_{m=1}^{M}$ approximating the posterior filtering density
$\pi_{n+1|n+1}$. The filtered state estimate is the posterior ensemble mean
\begin{equation}
    \widehat X_{t_{n+1}}
    =
    \frac{1}{M}
    \sum_{m=1}^{M}
    x_{n+1|n+1}^{(m)}.
    \label{eq:ensf_state_estimate}
\end{equation}

\section{Numerical Implementation}
\label{sec:numerical_implementation}

\subsection{Customized transmission network environment}
\label{sec:numerical_grid}

The numerical experiments use a customized transmission network simulator inspired by
Grid2Op. Its topology is derived from the
\texttt{l2rpn\_icaps\_2021\_large} environment, which contains 36 substations,
59 transmission lines, 37 load devices, and 22 generators. In Grid2Op, a
substation may contain multiple busbars, and different busbar assignments can
change the electrical topology within that substation. The customized simulator
does not retain this internal busbar topology. Each substation is instead treated
as a single electrical node.

The transmission line connectivity and the locations of loads and generators are
retained from the customized network snapshot. Active power flows are computed
using the DC power flow approximation described in
Section~\ref{sec:grid_control_problem}. The line parameters and active power limits
are numerical parameters of the customized DC model and should not be interpreted
as a reconstruction of the complete AC network represented by the original
Grid2Op environment. The simplified representation preserves the spatial network
structure needed for transmission constraints, line outages, storage, and generator
control while omitting the full busbar topology action space.

\begin{figure}[p]
    \centering
\includegraphics[width=0.88\linewidth]    {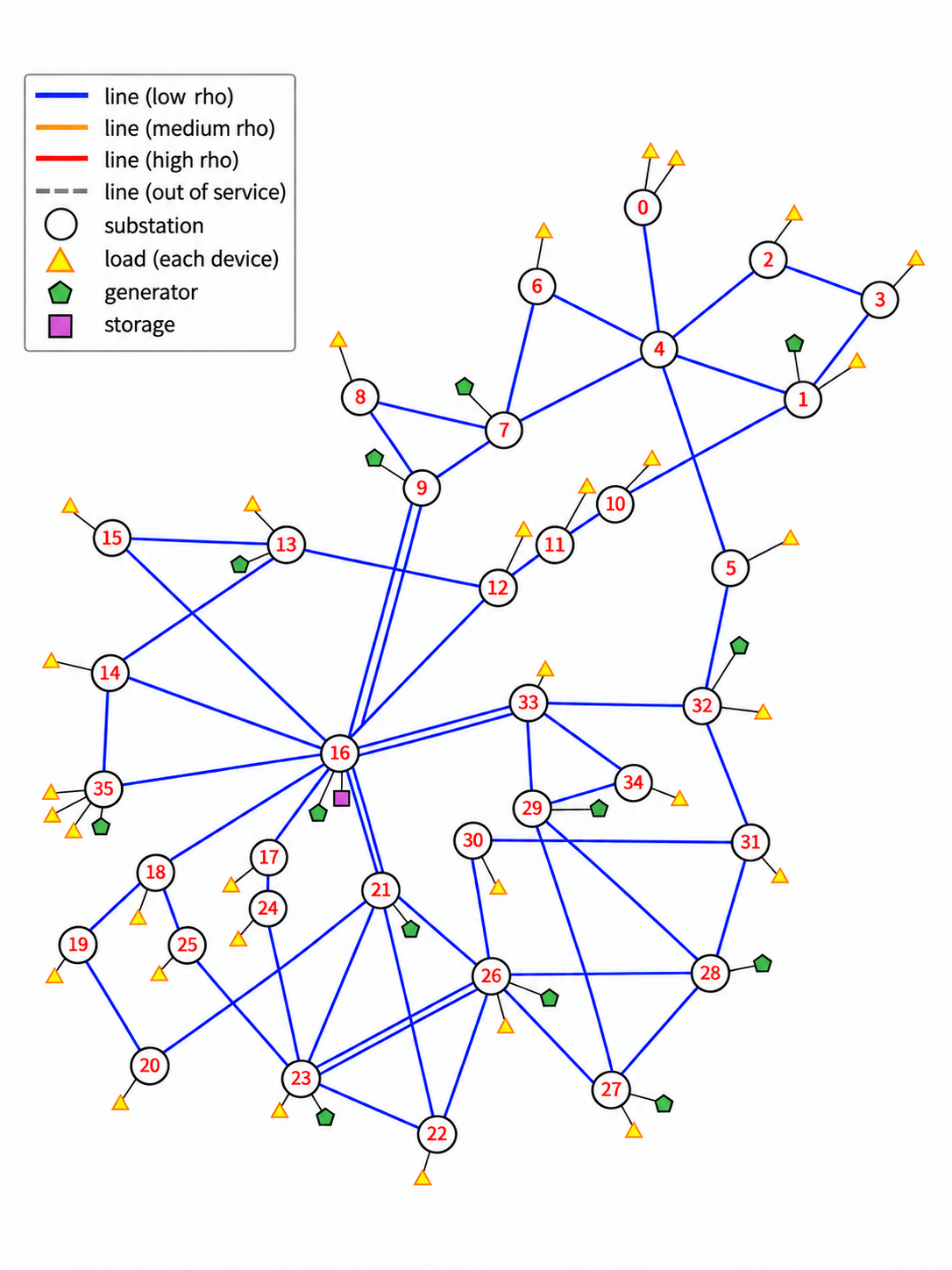}
\vspace{-0.8cm}
    \caption{
    Customized transmission network used in the numerical experiments.
    The network contains 36 substations, 59 transmission lines, 37 load devices,
    22 generators, and one energy storage device. Internal busbar topology is not
    represented in the customized simulator. The substation positions are schematic
    and do not represent geographical locations.
    }
    \label{fig:power_grid_topology}
\end{figure}

Figure~\ref{fig:power_grid_topology} shows the resulting network and the locations
of the load, generation, and storage devices.

\subsection{Network evolution, line protection, and storage}
\label{sec:numerical_network_dynamics}

The simulator uses a physical time step of 5 min. At each time $t_n$, the true
load and other exogenous injections are assigned to the physical network before
the controller selects an action using the information available at that time.
Generator control and any admissible line reconnection request are applied before
storage balancing. The resulting active power injections are used to solve the DC
power flow problem. If line protection changes the network topology, the power
flow is recomputed on the updated network.

Transmission line protection is modeled using an overload timer. For each online
line, the timer increases while its loading ratio remains above the prescribed
threshold and is reset when the loading ratio returns below that threshold. A line
is disconnected when the sustained overload duration reaches the specified
protection time. The disconnected line then enters a cooldown period during which
reconnection is not permitted. After the cooldown expires, the line becomes
eligible for reconnection by the controller. The overload threshold, protection
time, and cooldown duration are numerical parameters.

A single storage device is connected to substation 16. After generator control is
applied, storage absorbs a remaining generation surplus or supplies a remaining
generation deficit subject to its power and energy limits. Charging and
discharging update the stored energy according to the storage model in
Eq.~\eqref{eq:storage_energy}. The storage capacity, initial energy, power limits,
and efficiency are treated as experiment parameters.

\subsection{Load information and EnSF implementation}
\label{sec:numerical_ensf}

The physical load and the load information supplied to the controller are
maintained as separate quantities. Let
$L_n^{\mathrm{true}}
=
(L_{1,n}^{\mathrm{true}},\ldots,L_{N_d,n}^{\mathrm{true}})$
denote the load applied to the physical network, and let $L_n^f$ denote the
corresponding forward load information. Here $N_d=37$ is the number of load
devices. Loads at the device level are aggregated at the 36 substations for the power
flow calculation.

The EnSF implementation follows the filtering formulation in
Section~\ref{sec:ensf}. Since the load forecast is available as an algebraic
trajectory rather than through an explicit state transition model, each load
component is represented by the normalized deviation
\begin{equation}
    x_{i,n}
    =
    \frac{L_{i,n}}{L_{i,\mathrm{nom}}}-1,
    \label{eq:normalized_load}
\end{equation}
where $L_{i,\mathrm{nom}}>0$ is the nominal value of load component $i$.
The normalized forward and true load states are denoted by $x_n^f$ and
$x_n^{\mathrm{true}}$, respectively.

At the initial assimilation time, the prior ensemble is centered on the
normalized forward load,
\begin{equation}
    x_{0|-1}^{(m)}
    =
    x_0^f+\xi_0^{(m)},
    \qquad
    m=1,\ldots,M.
    \label{eq:ensf_initial_prior}
\end{equation}
At later assimilation times, the prior ensemble is obtained by translating the
previous posterior ensemble according to the change in the forward load,
\begin{equation}
    x_{n|n-1}^{(m)}
    =
    x_{n-1|n-1}^{(m)}
    +
    \left(
        x_n^f-x_{n-1}^f
    \right)
    +
    \xi_n^{(m)}.
    \label{eq:ensf_algebraic_prediction}
\end{equation}
The process perturbations are sampled as
$\xi_n^{(m)}\sim\mathcal N(0,\sigma_{\mathrm{proc}}^2 I)$, where $I$ is the
identity matrix of the corresponding state dimension. The prior ensemble follows
the same forward load trajectory used by the Forward case while retaining information
from previous filtering updates.

Noisy nonlinear observations of the true normalized load are generated
componentwise as
\begin{equation}
    y_{i,n}
    =
    \arctan
    \left(
        x_{i,n}^{\mathrm{true}}
    \right)
    +
    \eta_{i,n},
    \qquad
    \eta_{i,n}
    \sim
    \mathcal N
    \left(
        0,
        \sigma_{\mathrm{obs}}^2
    \right),
    \label{eq:ensf_arctan_observation}
\end{equation}
where $\sigma_{\mathrm{obs}}$ is the observation noise standard deviation.
For this observation model, the likelihood score used in the EnSF update is
\begin{equation}
    \nabla_x\log\ell_n(x)
    =
    -
    \frac{
        \arctan(x)-y_n
    }{
        \sigma_{\mathrm{obs}}^2(1+x^2)
    },
    \label{eq:ensf_arctan_score}
\end{equation}
where the operations are applied componentwise. The true load enters the filter
only through the noisy observation in
Eq.~\eqref{eq:ensf_arctan_observation} and is not supplied directly to the
filtering update.

The numerical implementation uses $M=50$ ensemble members,
$\sigma_{\mathrm{obs}}=0.1$, and $\sigma_{\mathrm{proc}}=0.03$. The score
approximation uses a mini-batch size of $J=1$. The EnSF load estimate is obtained
from the posterior ensemble mean and transformed back to physical units,
\begin{equation}
    \widehat L_{i,n}^{\mathrm{EnSF}}
    =
    L_{i,\mathrm{nom}}
    \left(
        1+
        \frac{1}{M}
        \sum_{m=1}^{M}
        x_{i,n|n}^{(m)}
    \right).
    \label{eq:ensf_load_reconstruction}
\end{equation}

Three load information modes are considered,
\begin{equation}
    L_n^{\mathrm{vis}} 
    =
    \begin{cases}
        L_n^f,
        & \text{Forward},\\[2mm]
        \widehat L_n^{\mathrm{EnSF}},
        & \text{EnSF},\\[2mm]
        L_n^{\mathrm{true}},
        & \text{Truth}.
    \end{cases}
    \label{eq:load_information_modes}
\end{equation}
The physical network evolves under $L_n^{\mathrm{true}}$ in all three cases.
Only the load information supplied to the controller is changed.

\section{Numerical Examples}
\label{sec:numerical_examples}

Two numerical examples are considered. The first uses a fixed heuristic controller
to isolate the effect of load correction on grid operation. The second considers a
richer sequential control setting in which data assimilation is combined with a
reinforcement learning controller.

Load estimation errors are summarized using the aggregation convention of
each experiment. For an evaluated sequence containing $T$ samples and $d_L$
load components, define the instantaneous spatial RMSE as
\begin{equation}
    e_n
    =
    \left[
        \frac{1}{d_L}
        \sum_{i=1}^{d_L}
        \left(
            L_{i,n}^{\mathrm{vis}}
            -
            L_{i,n}^{\mathrm{true}}
        \right)^2
    \right]^{1/2}.
    \label{eq:instantaneous_load_rmse}
\end{equation}
Experiment~1 reports the pooled device RMSE
$\left(T^{-1}\sum_{n=1}^{T}e_n^2\right)^{1/2}$ and the total-load RMSE
\begin{equation}
    \left[
        \frac{1}{T}
        \sum_{n=1}^{T}
        \left(
            \sum_{i=1}^{d_L}
            \left(
                L_{i,n}^{\mathrm{vis}}
                -
                L_{i,n}^{\mathrm{true}}
            \right)
        \right)^2
    \right]^{1/2}.
    \label{eq:total_load_rmse}
\end{equation}
Experiment~2 reports the time-averaged substation RMSE,
$T^{-1}\sum_{n=1}^{T}e_n$.
Scenario summaries are averaged with equal weight within each reported
comparison. Unless a common plotting window is specified, trajectory-based
error summaries use the samples available before the corresponding rollout
terminates or reaches its horizon. When stopping times differ, these summaries
describe the realized trajectories rather than estimation accuracy over a
common time window.

\subsection{EnSF-Assisted Heuristic Control}
\label{sec:heuristic_results}

The heuristic control example is designed to examine whether correcting uncertain
load information with EnSF produces a measurable operational benefit when the
controller itself is held fixed. Forward, EnSF-corrected, and Truth load information
are supplied to the same generator feedback rule under paired physical conditions.
The Truth case provides an oracle information benchmark, while a Do Nothing case
with fixed generation is included as an uncontrolled reference. This design separates
the effect of load estimation from changes in the control strategy.

The forward trajectory retains the individual profiles of the 37 load devices. The
physical load is generated by applying a common time-varying stochastic multiplier,
\begin{equation}
    L_{i,n}^{\mathrm{true}}
    =
    S_n L_{i,n}^{f},
    \qquad
    i=1,\ldots,N_d,
    \label{eq:heuristic_truth_load}
\end{equation}
where $S_n \in (0.05,2.5)$ is a bounded smooth stochastic process and $N_d=37$. The same physical load
trajectory is used by all control cases within each paired scenario. The filtering
problem also includes the nonlinear noisy observations and process perturbations
described in Section~\ref{sec:numerical_ensf}. Consequently, the EnSF estimate is
not expected to reproduce the physical load exactly. The relevant question here is
whether the reduction in load estimation error is sufficient to improve subsequent
grid operation.

Five generators are continuously controlled between 0 and 100 MW. The remaining
17 generators are held at 40 MW, giving a fixed generation contribution of 680 MW.
For each controlled case, the visible total load is used to determine the required
controllable generation after subtracting this fixed contribution. The resulting
generation request is restricted to the feasible range and allocated equally among
the five controllable generators. Thus, the Forward, EnSF, and Truth cases differ
only in the load information entering the same feedback rule.

The storage device has an energy capacity of 4000 MWh and is initialized at
2000 MWh. Its maximum charging and discharging powers are both 1000 MW, with
storage efficiency 0.95. Storage is not explicitly managed by the heuristic
controller. Instead, it passively absorbs the remaining generation surplus or
supplies the remaining deficit after generator control. Line protection is evaluated
with loading thresholds of $1.00$, $1.05$, and $1.10$. A sustained overload for
30 min disconnects a line, followed by a 240 min cooldown period. Five stochastic
load realizations are considered at each threshold, giving 15 paired scenarios for
each information mode. A trajectory terminates if the storage boundary prevents
power balancing, if the remaining generation and load mismatch exceeds the admissible
tolerance, or if the active network no longer admits a valid DC power flow.

The primary operational quantity in this example is survival time. Load information
quality is evaluated using both root-mean-square error (RMSE) at the device level and for total load. Table~\ref{tab:heuristic_control_summary} summarizes the results over the
15 paired scenarios.
\begin{table}[h!]
    \centering
    \caption{
    Heuristic control performance over 15 paired scenarios. Survival values are
    reported as mean $\pm$ standard deviation. RMSE is not defined for the
    Do Nothing case because no load estimate is used by that controller.
    }
    \label{tab:heuristic_control_summary}
    \begin{tabular*}{\linewidth}{@{\extracolsep{\fill}}lccc@{}}
        \toprule
        Method
        & Survival
        & \makecell{Device RMSE\\(MW)}
        & \makecell{Total RMSE\\(MW)} \\
        \midrule
        Do Nothing
        & $316.3 \pm 97.2$
        & -- & -- \\
        Forward
        & $1336.5 \pm 645.7$
        & 7.10 & 151.33 \\
        EnSF
        & $2260.7 \pm 1395.4$
        & 5.67 & 118.91 \\
        Truth
        & $2540.1 \pm 897.6$
        & 0.00 & 0.00 \\
        \bottomrule
    \end{tabular*}
\end{table}

EnSF reduces the mean device level RMSE from 7.10 MW to 5.67 MW relative to the
Forward case, a reduction of approximately 20.2\%. The total load RMSE decreases
from 151.33 MW to 118.91 MW, corresponding to a reduction of approximately 21.4\%.
This improvement in load information produces a considerably larger change in the
operational outcome. Mean survival increases from 1336.5 steps with Forward
information to 2260.7 steps with EnSF correction, an increase of approximately
69.2\%. The Truth case gives the highest mean survival at 2540.1 steps.

\begin{figure}[!ht]
    \centering
\includegraphics[width=0.7\linewidth]    {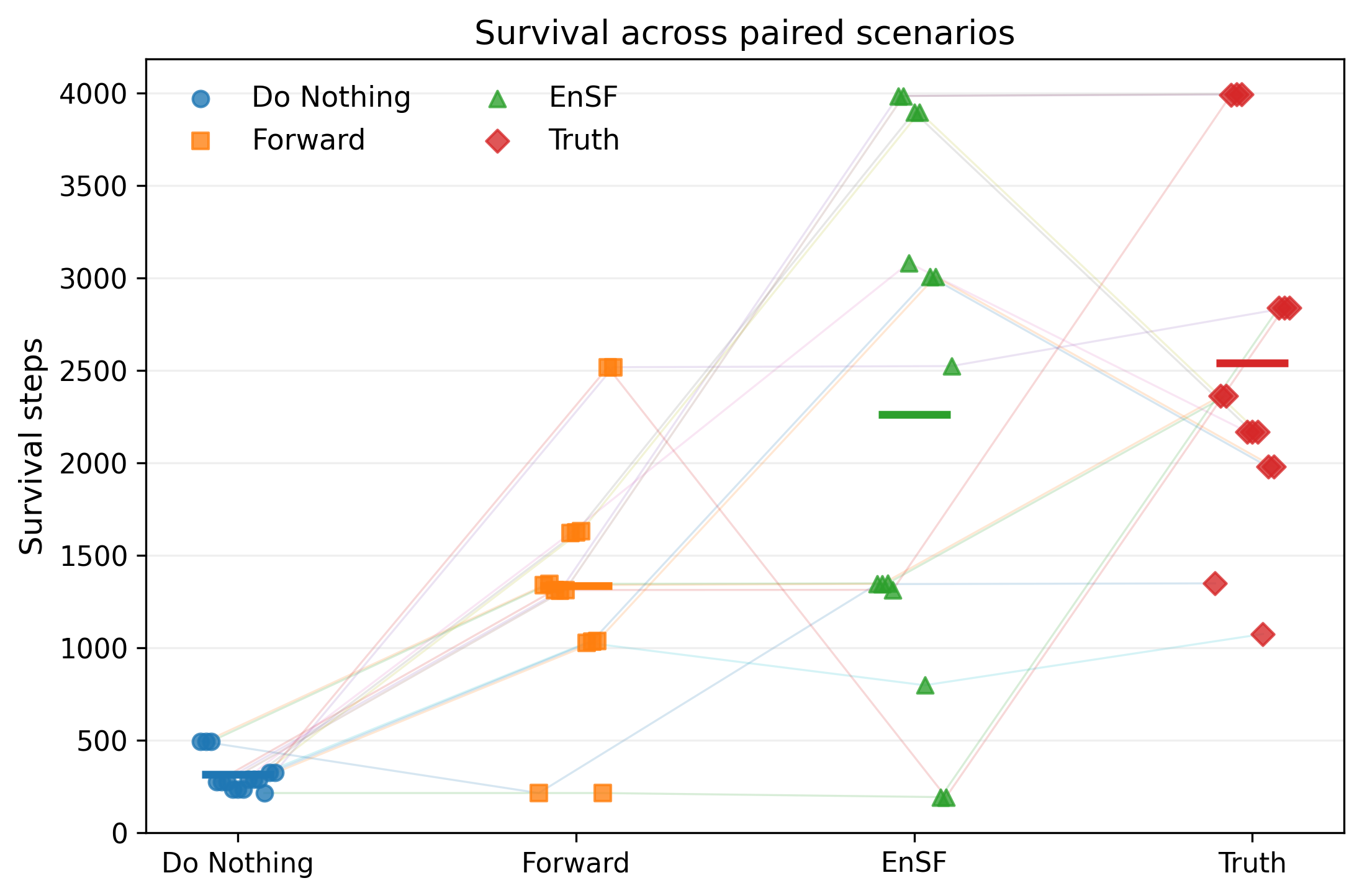}
    \caption{
    Survival across the 15 paired heuristic control scenarios. Each faint line
    connects the four cases evaluated under the same stochastic load realization
    and line protection threshold. The larger horizontal markers indicate the
    mean survival for each method.
    }
    \label{fig:heuristic_survival_paired}
\end{figure}

Figure~\ref{fig:heuristic_survival_paired} shows the individual paired outcomes.
Each connecting line corresponds to the same stochastic load realization and line
protection threshold under the four control cases. The large variation across
scenarios reflects the nonlinear interaction among load variation, transmission
constraints, line protection, and storage. EnSF generally shifts the controlled
trajectories toward longer survival relative to Forward information.

Individual EnSF trajectories can occasionally survive longer than the corresponding
Truth trajectory. This does not indicate that the EnSF estimate contains more accurate
information than the oracle case. The heuristic controller does not explicitly manage
storage state of charge, and storage passively balances the mismatch left after
generator dispatch. Differences in the load information can change earlier
charging and discharging behavior and, in some realizations, leave different amounts
of storage headroom before a later generation surplus. Truth should consequently be
interpreted as an information benchmark rather than an optimal control benchmark in
this example.

\begin{figure}[!ht]
    \centering
    \includegraphics[width=0.96\linewidth]
    {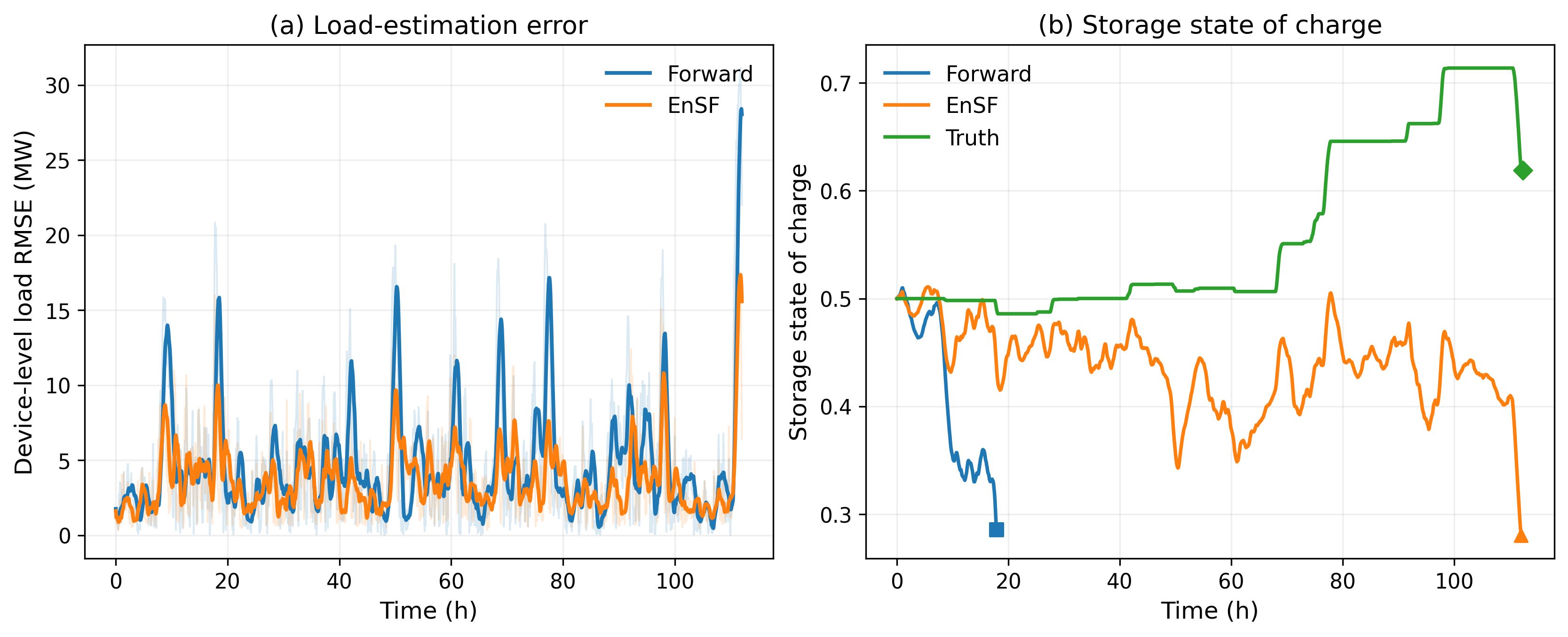}
    \caption{
    Representative heuristic control trajectory at line loading threshold $1.00$.
    The left figure compares load RMSE at the device level for the Forward and EnSF
    information channels, with thin curves showing instantaneous RMSE and thicker
    curves showing 60-min moving averages. The right figure shows the corresponding
    storage state of charge trajectories under Forward, EnSF, and Truth information.
    }
    \label{fig:heuristic_filter_storage}
\end{figure}

Figure~\ref{fig:heuristic_filter_storage} shows the connection between filtering
and control for one representative paired scenario with line loading threshold $1.00$.
EnSF reduces the load estimation error relative to the forward model over much of the
trajectory. Because the same feedback rule acts on the corrected load information,
this difference propagates into generator dispatch and storage use.

The EnSF gain is not confined to the estimation metric. Despite model discrepancy,
nonlinear noisy observations, and process uncertainty, the corrected load information
substantially extends survival under the same heuristic rule. The occasional reversal
between EnSF and Truth also exposes a limitation of the controller: it reacts to current
estimated demand without explicitly managing storage state or future consequences.

\FloatBarrier

\subsection{Data-Assimilation-Assisted Reinforcement Learning}
\label{sec:ppo_results}

The second example adds sequential discrete control. Generator dispatch, line
reconnection, renewable variability, storage, and load uncertainty now interact over
time, and the controller acts on a structured grid observation. The main questions are
whether better load information still helps after the control law has been learned and
whether EnSF can limit performance loss when a frozen policy is evaluated under stronger
load uncertainty than it encountered during training.

The reinforcement learning configuration contains 10 dispatchable generators and
12 nondispatchable renewable-like generators. The dispatchable generators are located at
substations
\[
    1,\ 7,\ 9,\ 13,\ 16,\ 23,\ 27,\ 29,\ 32,\ 35.
\]
Each renewable-like generator produces between 0 and 60 MW. The stochastic output is
correlated with persistence coefficient 0.97 and mean operating fraction 0.5.
The common and generator-specific perturbations have standard deviations 0.020 and
0.025, respectively, and the initial perturbation standard deviation is 0.080. Realized
renewable generation is visible to the controller in all load information modes. For a
paired comparison, the same physical load and renewable trajectories are used under
Forward, EnSF, and Truth information.

The storage device has a capacity of 2000 MWh and is initialized at 500 MWh. Its maximum
charging and discharging powers are 1000 MW, and the storage efficiency parameter is
0.95. A line is considered overloaded when its loading ratio exceeds one. An overload
sustained for 30 min causes the line to trip, after which a 240 min cooldown is imposed
before reconnection becomes available. A rollout terminates if the active network becomes
singular or disconnected, the residual generation and load mismatch exceeds 1 MW, or the
storage energy reaches an enforced boundary. Each evaluation rollout has a maximum
length of 3000 time steps.

For reinforcement learning, each dispatchable generator is restricted to the five
active power levels
\begin{equation}
    \mathcal P_g
    =
    \left\{
        0,\ 20,\ 40,\ 60,\ 80
    \right\}
    \ \mathrm{MW}.
    \label{eq:ppo_generator_levels}
\end{equation}
The action at time $t_n$ is
\begin{equation}
    A_n
    =
    \left(
        a_{1,n}^{g},
        \ldots,
        a_{10,n}^{g},
        a_n^{\mathrm{rec}}
    \right),
    \label{eq:ppo_action}
\end{equation}
where $a_{k,n}^{g}\in\{0,1,2,3,4\}$ selects one of the generator levels in
Eq.~\eqref{eq:ppo_generator_levels}. The reconnection component
$a_n^{\mathrm{rec}}\in\{0,1,\ldots,59\}$ selects no reconnection when
$a_n^{\mathrm{rec}}=0$ and requests reconnection of line
$a_n^{\mathrm{rec}}-1$ otherwise. A reconnection is applied only when the selected line
is offline and its cooldown has expired.

The PPO observation is the 187-dimensional vector
\begin{equation}
    O_n
    =
    \left[
        O_n^{\mathrm{load}},
        O_n^{\mathrm{ren}},
        O_n^{\mathrm{storage}},
        O_n^{\mathrm{summary}},
        O_n^{\mathrm{line}},
        O_n^{\mathrm{status}},
        O_n^{\mathrm{gen}}
    \right].
    \label{eq:ppo_observation}
\end{equation}
The load block $O_n^{\mathrm{load}}\in\mathbb R^{37}$ contains the normalized total
visible load and 36 spatial substation load values.
The renewable block $O_n^{\mathrm{ren}}\in\mathbb R^{12}$ contains the realized
renewable-like generation. The storage block has two components, state of charge and
normalized storage power. The eight summary variables contain the three largest line
loading ratios, the fractions of overloaded and offline lines, the previous maximum line
loading, and sine and cosine encodings of time of day. The observation also contains
59 line loading ratios, 59 line control states, and the normalized operating levels of
the 10 dispatchable generators. When the load information mode changes, only
$O_n^{\mathrm{load}}$ is replaced. The remaining observation components are generated
from the same physical realization.

Before PPO is evaluated, several heuristic reference controllers are considered in the
same environment. The Forward, EnSF, and Truth heuristics greedily select the discrete
generator levels to reduce the instantaneous mismatch between controllable generation
and the visible net load. Storage then absorbs the remaining surplus or supplies the
remaining deficit. These heuristics do not otherwise regulate the storage state of charge.

A storage-aware Truth reference is included to distinguish accurate load information from
the limitation of dispatch based only on instantaneous mismatch. Let $\overline E=2000$ MWh denote the
storage capacity and let $E_n$ denote the stored energy. The storage recovery power is
defined by
\begin{equation}
    P_n^{\mathrm{rec}}
    =
    \operatorname{clip}
    \left(
        \frac{0.5\overline E-E_n}{12~\mathrm{h}},
        -100,\ 100
    \right)
    \ \mathrm{MW}.
    \label{eq:storage_recovery}
\end{equation}
The controllable generation target is then
\begin{equation}
    P_n^{\mathrm{target}}
    =
    P_n^{\mathrm{true}}
    -
    P_n^{\mathrm{ren}}
    +
    P_n^{\mathrm{rec}},
    \label{eq:storage_aware_target}
\end{equation}
where $P_n^{\mathrm{true}}$ is the total physical load and
$P_n^{\mathrm{ren}}$ is the realized renewable generation. The correction in
Eq.~\eqref{eq:storage_recovery} drives the storage gradually toward 50\% state of charge
and is limited to 100 MW. This reference addresses the storage drift that can occur when
generation is chosen only to match the current load.

Table~\ref{tab:rl_heuristic_summary} summarizes the heuristic reference results over five
paired physical scenarios.

\begin{table}[!ht]
    \centering
    \caption{
    Heuristic reference controllers in the reinforcement learning environment. The reference controllers are scored using the same per-step reward used for PPO in Eq.~\eqref{eq:ppo_reward}.
    Survival and return are reported as mean $\pm$ standard deviation over five paired
    physical scenarios.
    }
    \label{tab:rl_heuristic_summary}
    \begin{tabular*}{\linewidth}{@{\extracolsep{\fill}}lcccc@{}}
        \toprule
        Controller
        & Survival
        & Return
        & \makecell{Warning\\fraction}
        & \makecell{Overload\\fraction} \\
        \midrule
        Do Nothing
        & $32.4 \pm 19.0$
        & $-49.9 \pm 16.4$
        & 0.8684
        & 0.4611 \\
        Forward
        & $679.2 \pm 541.2$
        & $350.3 \pm 376.9$
        & 0.1247
        & 0.0157 \\
        EnSF
        & $1627.0 \pm 627.5$
        & $575.4 \pm 178.8$
        & 0.0612
        & 0.0066 \\
        Truth
        & $2106.2 \pm 1006.6$
        & $838.5 \pm 569.7$
        & 0.0327
        & 0.0003 \\
        Truth + Storage
        & $3000.0 \pm 0.0$
        & $2952.3 \pm 13.0$
        & 0.0250
        & 0.0002 \\
        \bottomrule
    \end{tabular*}
\end{table}

The heuristic references preserve the same ordering in load information quality. Mean
survival rises from 679.2 steps with Forward information to 1627.0 with EnSF and 2106.2
with Truth, while the mean substation load RMSE falls from 5.35 MW to 4.00 MW after
EnSF correction. The storage-aware Truth controller reaches the full 3000-step horizon
in all five scenarios. Accurate load information alone is not enough to manage
storage when dispatch is purely myopic. Direct attention to state of charge removes the
drift seen with the simpler rule.

For PPO training, let $P_n^{\mathrm{mis}}$ denote the generation and load mismatch before
storage balancing, let $\delta_n$ denote the residual mismatch after storage balancing,
and let $P_n^s$ denote storage power. The state of charge is
\[
    q_n=\frac{E_n}{\overline E}.
\]
Let $N_n^{\mathrm{warn}}$ be the number of online lines with
$\rho_{\ell,n}>0.8$, and let $N_n^{\mathrm{over}}$ be the number with
$\rho_{\ell,n}>1$. The overload severity is
\begin{equation}
    V_n
    =
    \sum_{\ell\in\mathcal E_n^{\mathrm{on}}}
    \left(
        \rho_{\ell,n}-1
    \right)_+,
    \label{eq:overload_severity}
\end{equation}
where $\mathcal E_n^{\mathrm{on}}$ is the set of online lines and
$(z)_+=\max(z,0)$. Let $N_n^{\mathrm{trip}}$ denote the number of newly tripped
lines, let $I_n^{\mathrm{rec}}$ indicate a successful reconnection, and let
$I_n^{\mathrm{term}}$ indicate episode termination. With
$P_{\mathrm{nom}}^d$ denoting total nominal load and $P_{\max}^s=1000$ MW denoting
the storage power scale, the PPO reward is
\begin{equation}
\begin{aligned}
    r_n ={}&
    1
    -
    \frac{|P_n^{\mathrm{mis}}|}{P_{\mathrm{nom}}^d}
    -
    0.10\frac{|\delta_n|}{P_{\mathrm{nom}}^d}
    -
    0.05\frac{|P_n^s|}{P_{\max}^s} \\[1mm]
    &-
    10\left[
        (0.4-q_n)_+^2
        +
        (q_n-0.6)_+^2
    \right]
    -
    0.01N_n^{\mathrm{warn}}
    -
    0.05N_n^{\mathrm{over}} \\[1mm]
    &-
    0.25V_n
    -
    0.50N_n^{\mathrm{trip}}
    -
    0.01I_n^{\mathrm{rec}}
    -
    50I_n^{\mathrm{term}}.
    \label{eq:ppo_reward}
\end{aligned}
\end{equation}
The reward penalizes dispatch mismatch, storage use, excursions of state of charge outside
the interval $[0.4,0.6]$, transmission stress, line trips, line reconnection and terminal failure. The SOC
term gives the learned controller an explicit incentive to retain storage flexibility over
the rollout rather than using storage only as an instantaneous balancing device.

\begin{figure}[!ht]
    \centering
    \begin{subfigure}[t]{\linewidth}
        \centering
        \includegraphics[width=0.7\linewidth]
        {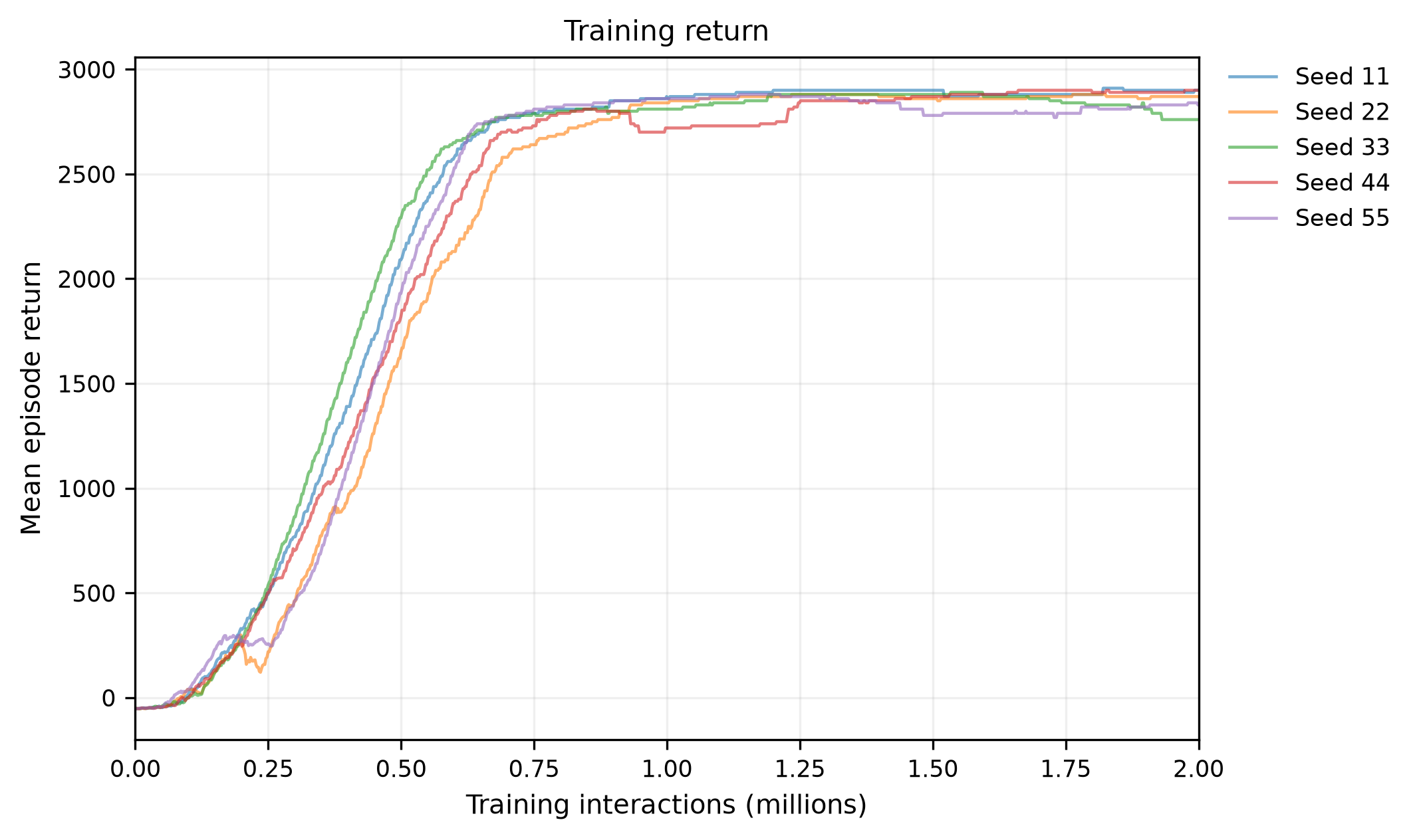}
        \caption{Mean training return.}
        \label{fig:ppo_training_return}
    \end{subfigure}

    \vspace{0.08in}

    \begin{subfigure}[t]{0.76\linewidth}
        \centering
        \includegraphics[width=0.95\linewidth]
        {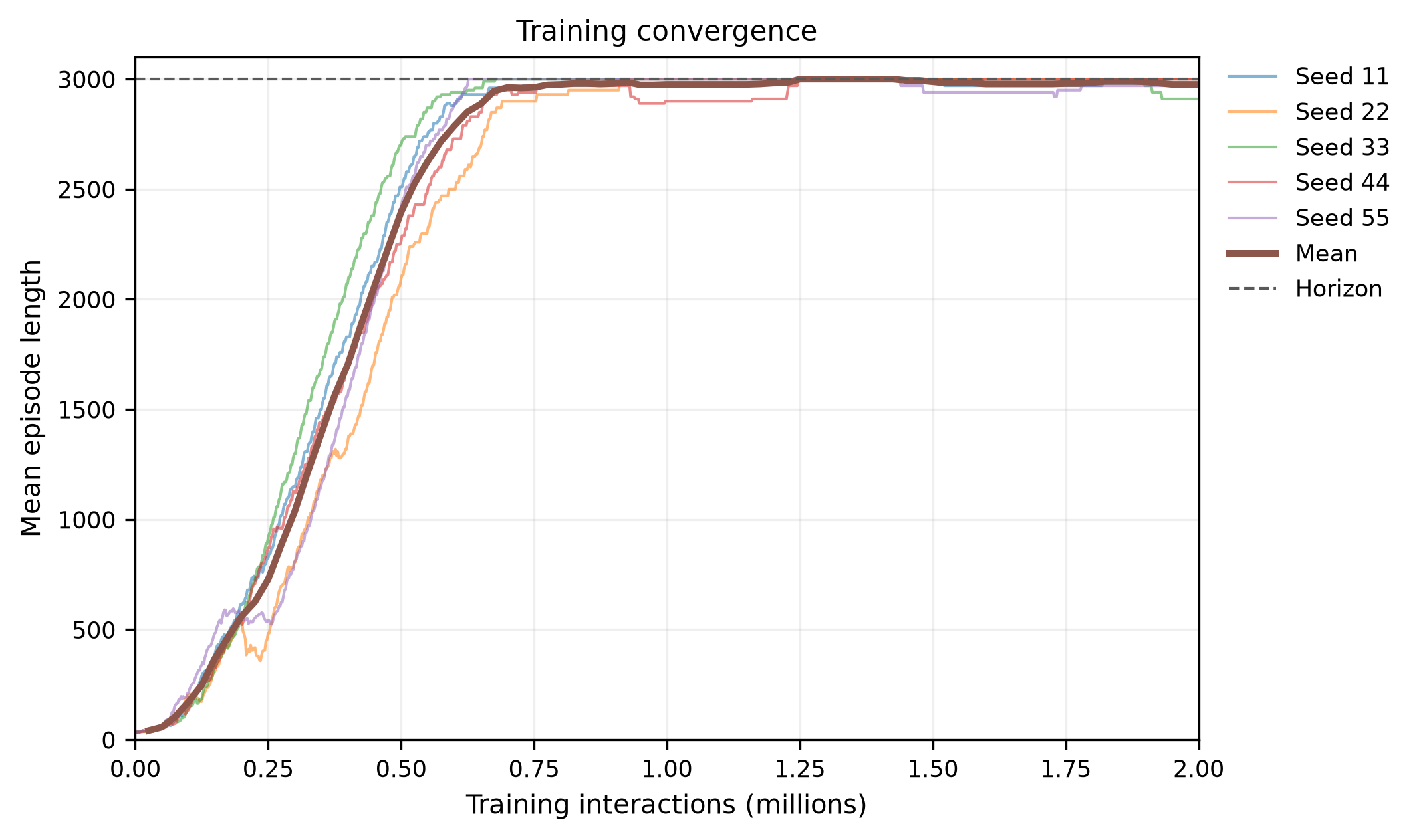}
        \caption{Mean episode length.}
        \label{fig:ppo_training_length}
    \end{subfigure}
    \caption{
    PPO training over two million environment interactions for five independent
    training runs. The horizontal line in the episode length panel marks the
    3000-step rollout horizon.
    }
    \label{fig:ppo_training}
\end{figure}

PPO \cite{schulman2017} is implemented with Stable-Baselines3
\cite{raffin2021}. The multilayer perceptron policy uses a learning rate of
$3\times10^{-4}$, rollout length 2048, and mini-batch size 64, with the remaining PPO
parameters set to the Stable-Baselines3 defaults. Five policies are trained independently
for $2\times10^6$ environment interactions under Truth load information. EnSF is not
used during training. Figure~\ref{fig:ppo_training} shows that the mean training return
and episode length approach a stable regime after approximately $8\times10^5$
interactions.

\FloatBarrier

After training, the policy parameters are frozen. Each policy is evaluated under Forward,
EnSF, and Truth load information using five paired physical scenarios. The resulting
evaluation contains 25 combinations of a trained policy and physical scenario for each information mode. All 75 rollouts
reach the 3000-step horizon, so survival does not distinguish the three information modes
under the nominal stochastic load distribution.

\begin{table}[!ht]
    \centering
    \caption{
    Frozen PPO evaluation under the nominal load distribution. Return values are
    reported over 25 combinations of a trained policy and physical scenario for each information mode. Load RMSE is
    determined by the five unique physical scenarios.
    }
    \label{tab:ppo_summary}
    {\setlength{\tabcolsep}{3pt}
    \begin{tabular*}{\linewidth}{@{\extracolsep{\fill}}lccccc@{}}
        \toprule
        Information
        & Return
        & \makecell{Time-averaged\\substation RMSE (MW)}
        & \makecell{Warning\\steps}
        & \makecell{Overload\\steps}
        & \makecell{Mean\\SOC} \\
        \midrule
        Forward
        & $2769.3 \pm 36.5$
        & 5.196
        & 57.60
        & 0.48
        & 0.497 \\
        EnSF
        & $2859.2 \pm 14.8$
        & 3.927
        & 44.72
        & 0.36
        & 0.488 \\
        Truth
        & $2913.1 \pm 29.8$
        & 0.000
        & 34.48
        & 0.24
        & 0.474 \\
        \bottomrule
    \end{tabular*}}
\end{table}
\FloatBarrier

EnSF increases the mean frozen policy return from 2769.3 to 2859.2, corresponding to
an improvement of approximately 3.25\% relative to Forward information. EnSF exceeds
Forward in all 25 paired comparisons. Truth information gives a mean
return of 2913.1 and exceeds EnSF in 24 of the 25 pairs. The mean substation load RMSE
decreases from 5.196 MW to 3.927 MW after EnSF correction, a reduction of approximately
24.4\%. Warning exposure also decreases from 57.60 to 44.72 steps per rollout. The
overload counts are already small in this nominal regime, so their differences are less
pronounced.

\begin{figure}[!ht]
    \centering
    \begin{subfigure}[t]{0.49\linewidth}
        \centering
        \includegraphics[width=\linewidth]
        {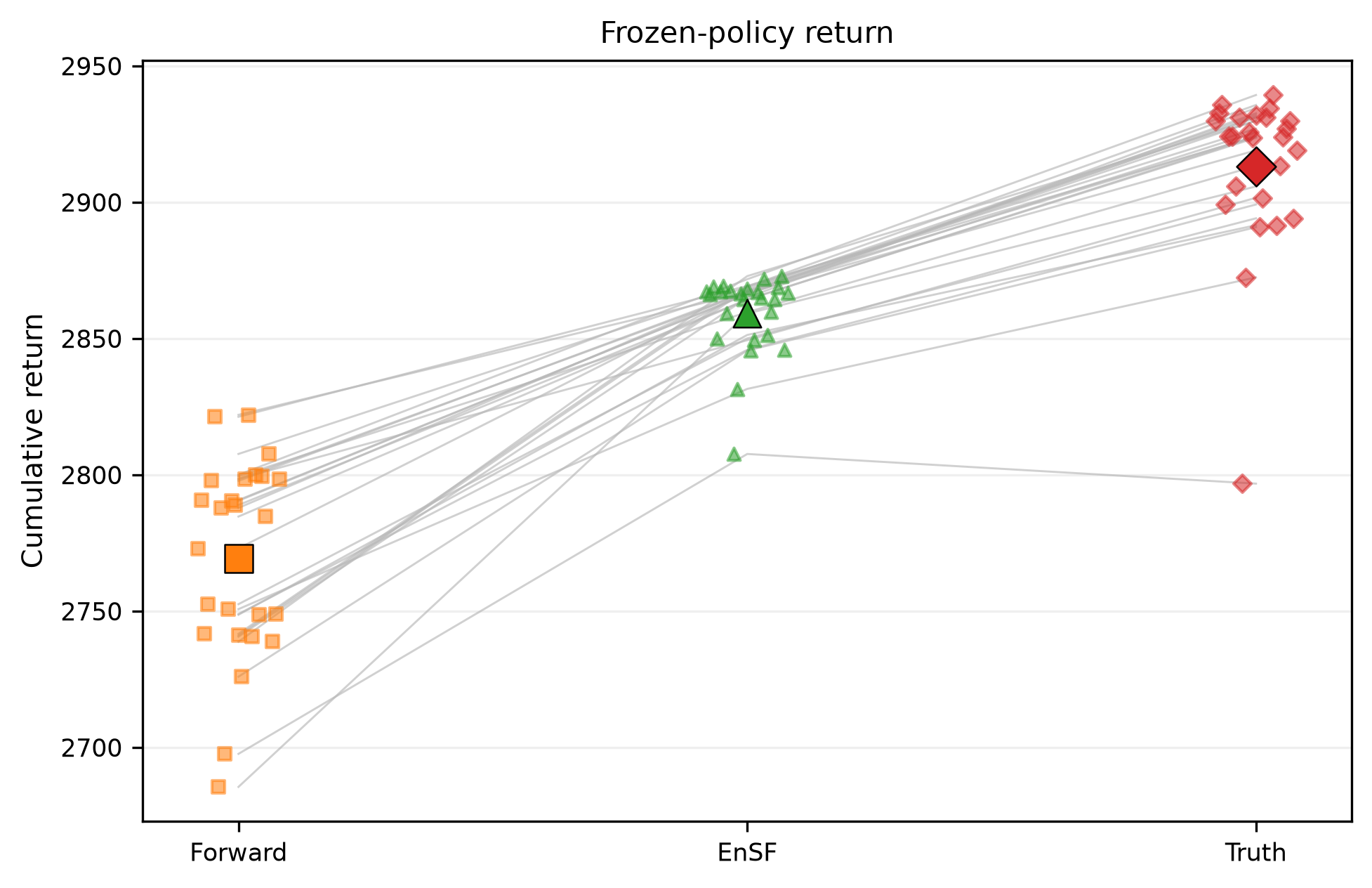}
        \caption{Frozen policy cumulative return.}
        \label{fig:ppo_return_paired}
    \end{subfigure}
    \hfill
    \begin{subfigure}[t]{0.49\linewidth}
        \centering
        \includegraphics[width=\linewidth]
        {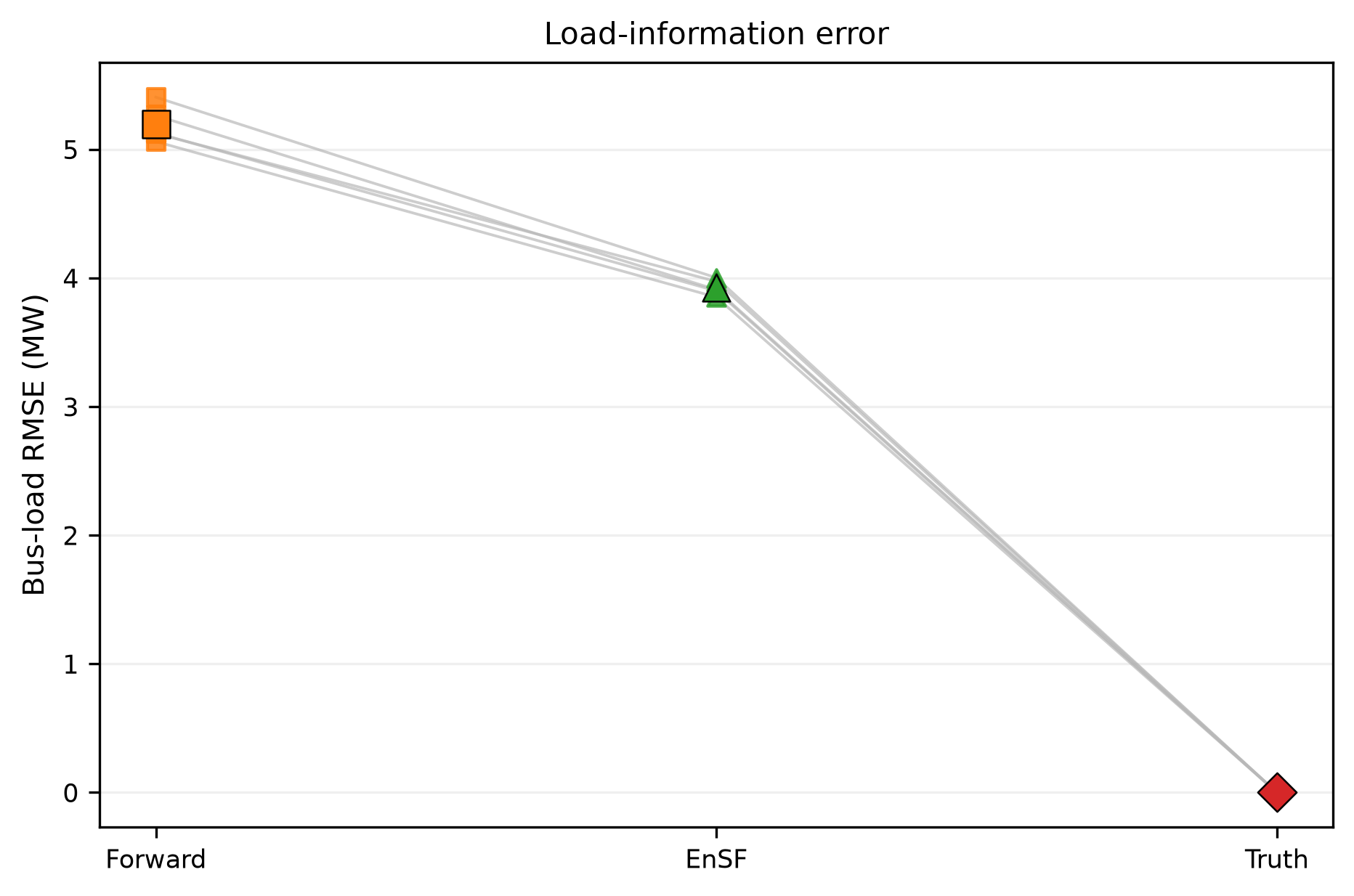}
        \caption{Substation load RMSE.}
        \label{fig:ppo_load_rmse}
    \end{subfigure}
    \caption{
    Paired frozen policy evaluation under the nominal load distribution. The return
    panel contains five trained policies evaluated on five physical scenarios. The
    load error panel represents the five distinct physical load trajectories.
    }
    \label{fig:ppo_paired_results}
\end{figure}
\FloatBarrier

Because survival reaches the evaluation ceiling in the nominal regime, a second
evaluation is performed under stronger load stochasticity without retraining the PPO
policies. Let $S_n$ denote the common stochastic load multiplier in the nominal load
model. The stressed multiplier is defined as
\begin{equation}
    S_n^{\mathrm{stress}}
    =
    1
    +
    \gamma
    \left(
        S_n-1
    \right),
    \qquad
    \gamma=1.40.
    \label{eq:ppo_load_stress}
\end{equation}
Thus, the departure of the common load multiplier from one is amplified by a factor of
1.4. The temporal realization and the existing spatial load structure are retained. The
forward trajectory, renewable process, grid parameters, storage parameters, EnSF
parameters, and PPO weights are unchanged. The value $\gamma=1.40$ is fixed on a separate calibration set using only Truth
information, before any Forward or EnSF stress results are evaluated.

The formal stress evaluation uses ten physical load realizations and the same five frozen
PPO policies, giving 50 rollouts for each information mode. Forward, EnSF, and Truth are
evaluated on identical stressed physical trajectories within each paired comparison.
Table~\ref{tab:ppo_stress_summary} summarizes the resulting performance.

\begin{table}[!ht]
    \centering
    \caption{
    Frozen PPO evaluation under increased load stochasticity. Survival and return are
    reported over 50 combinations of a trained policy and physical scenario for each information mode. Warning fraction
    is the fraction of survived steps containing at least one transmission warning.
    }
    \label{tab:ppo_stress_summary}
    {\setlength{\tabcolsep}{2pt}
    \begin{tabular*}{\linewidth}{@{\extracolsep{\fill}}lccccc@{}}
        \toprule
        Information
        & Survival
        & \makecell{Fraction reaching\\3000 steps}
        & Return
        & \makecell{Time-averaged substation\\load RMSE (MW)}
        & \makecell{Warning\\fraction} \\
        \midrule
        Forward
        & $2894.0 \pm 325.7$
        & 0.90
        & $2473.2 \pm 376.9$
        & 6.496
        & 0.0384 \\
        EnSF
        & $2963.3 \pm 252.8$
        & 0.96
        & $2750.8 \pm 254.2$
        & 4.820
        & 0.0319 \\
        Truth
        & $2964.1 \pm 253.6$
        & 0.98
        & $2808.7 \pm 270.2$
        & 0.000
        & 0.0257 \\
        \bottomrule
    \end{tabular*}}
\end{table}

Stronger load variability separates the information modes more clearly. Mean return
rises from 2473.2 with Forward information to 2750.8 with EnSF, an improvement of
approximately 11.2\%. The substation load RMSE falls by approximately 25.8\%, from
6.496 MW to 4.820 MW, and the mean warning fraction falls by approximately 16.9\%.
Forward reaches the full 3000-step horizon in 45 of 50 rollouts, compared with 48 of
50 for EnSF and 49 of 50 for Truth.

The five PPO policies share the same ten physical scenarios, so the 50 rollouts within
one information mode are not treated as 50 independent physical realizations. Averaging
the five policies within each physical scenario gives ten paired scenario level
comparisons. EnSF produces a higher mean return than Forward in nine of the ten scenarios.
A two-sided Wilcoxon signed-rank test applied to these ten paired scenario-level return differences gives \(p=0.0098\). The warning fraction is
lower with EnSF in all ten scenarios, with $p=0.0020$. Survival remains partially
censored by the 3000-step horizon and provides less statistical separation than
return and warning exposure.

\begin{figure}[!ht]
    \centering
    \begin{subfigure}[t]{0.49\linewidth}
        \centering
        \includegraphics[width=\linewidth]
        {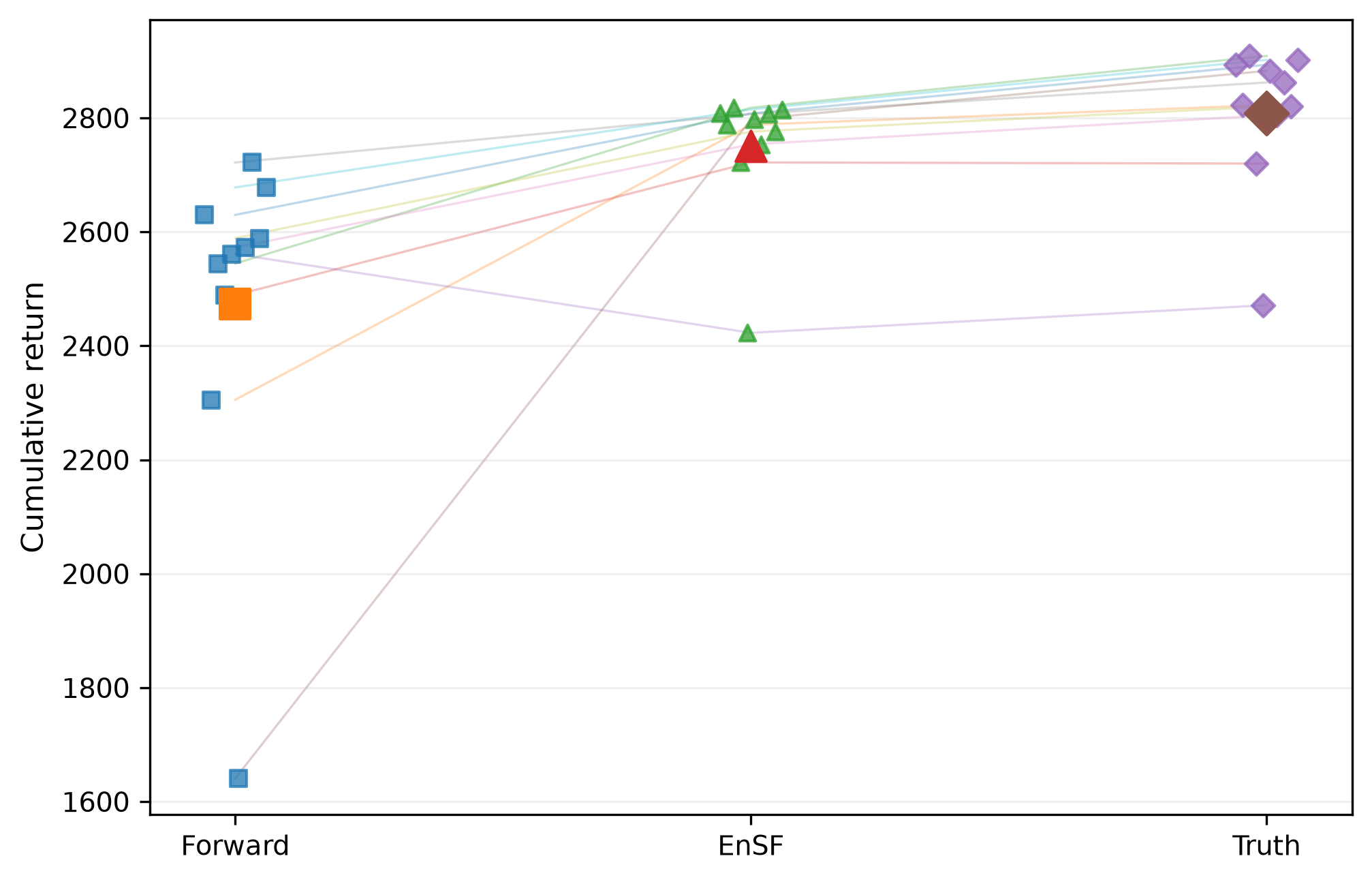}
        \caption{Paired cumulative return under increased load stochasticity.}
        \label{fig:ppo_stress_return}
    \end{subfigure}
    \hfill
    \begin{subfigure}[t]{0.49\linewidth}
        \centering
        \includegraphics[width=\linewidth]
        {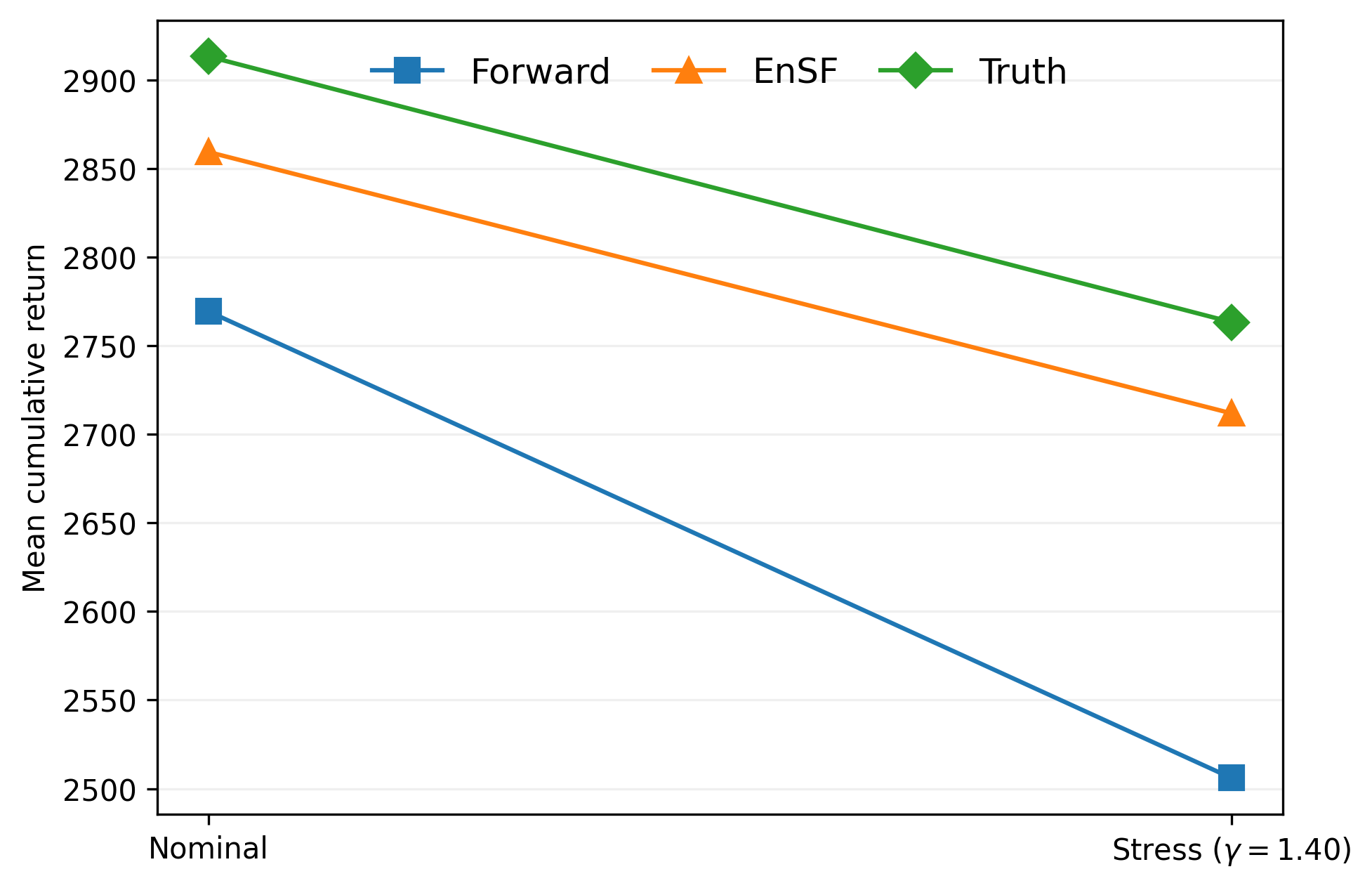}
        \caption{Matched nominal and stress evaluation.}
        \label{fig:ppo_nominal_stress_return}
    \end{subfigure}
    \caption{
    Frozen policy performance under the load distribution shift. In the left panel,
    each line represents one physical stress scenario after averaging over the five
    trained policies. The right panel compares mean return under nominal and stressed
    load conditions using only the physical scenarios shared by both evaluations.
    }
    \label{fig:ppo_stress_results}
\end{figure}

For the physical scenarios shared by the nominal and stress evaluations, the mean Forward
return decreases from 2769.3 to 2506.0, a reduction of approximately 9.6\%. The
corresponding decreases for EnSF and Truth are 5.2\% and 5.2\%, respectively. The
EnSF--Forward return gap therefore increases from approximately 90.0 in the nominal
regime to 205.8 under stress, while the Truth--EnSF gap remains nearly unchanged at
53.9 and 51.5. The larger deployment shift affects the Forward information channel much
more strongly than the EnSF and Truth channels.

Under nominal load variability, survival is saturated, but EnSF still improves return
and warning exposure. Increasing the load uncertainty without retraining the policy
widens the gap between Forward and EnSF, while EnSF remains close to the Truth
benchmark and fails early less often. Load correction matters more once the frozen PPO
policy is moved away from its nominal training distribution.

\FloatBarrier

\section{Conclusion}
\label{sec:conclusion}

We studied power grid control when the load information available to the controller is
uncertain and can be corrected by data assimilation. With the heuristic controller held
fixed, EnSF reduces device level load RMSE by about 20\% relative to the forward model
and increases mean survival by about 69\%. The frozen PPO policies also benefit from
EnSF under the nominal load distribution, where every paired comparison gives a higher
return than with Forward information. When load variability is increased beyond the
nominal training regime, the mean PPO return improves by about 11.2\% with EnSF and
the fraction of rollouts reaching the full horizon rises from 90\% to 96\%. In both
experiments, EnSF moves the controller toward the performance obtained with Truth
information.

The filtering results are consistent across two different control structures. EnSF does
more than reduce load estimation error: the corrected information changes generator
decisions, storage use, warning exposure, survival, and cumulative return. The stress
test also indicates that load estimate quality matters more when the operating
distribution shifts away from the one used for policy training. In
the present design, EnSF is applied after PPO training, so a trained policy can use
corrected load information without running the filter throughout millions of training
interactions. Truth remains an information benchmark rather than an optimal control
benchmark. Accurate state information cannot compensate for a controller that ignores
longer term quantities such as storage state of charge. The storage-aware heuristic and
the PPO reward in the second example show how state estimation and control design can
address different parts of the same problem.

The current numerical model is deliberately simplified. It uses DC power flow and does
not represent reactive power, voltage magnitude, frequency dynamics, or the full internal
busbar topology of the original Grid2Op environment. Load observations, spatial load
uncertainty, and renewable-like generation are synthetic, and renewable output is assumed
to be directly available to the controller. PPO uses a restricted set of discrete generator
levels and line reconnection actions. Storage responds through power balancing and the
reward rather than through an independent control action. The stress test changes only
the amplitude of load stochasticity, not renewable uncertainty, network outages, model
parameters, or observation availability.

Several extensions follow naturally from these limitations. Higher fidelity studies could
use AC or dynamic power system models, larger networks, measured load and renewable
data, richer topology actions, direct storage control, generator ramping, and multiple
sources of uncertainty. The computational cost of EnSF also needs to be studied at larger
online assimilation scales. Temperature and other weather variables could be incorporated
into the load model so that weather forecasts contribute to load prediction and subsequent
data assimilation. A separate direction is to compare the present modular design with
training schemes that place EnSF inside the reinforcement learning loop, as well as
belief-state, uncertainty-aware, or risk-sensitive control.

\section{Acknowledgements}

This research was supported by the U.S. Department of Energy (DOE), Office of Science, Office of Advanced Scientific Computing Research (ASCR), Applied Mathematics Program under contracts ERKJ443 and ERKJ388.
Feng Bao acknowledges support from the U.S. National Science Foundation (NSF) under Grants DMS-2142672 and DMS-2608431, and from the U.S. Department of Energy, Office of Science, Advanced Scientific Computing Research (ASCR) Applied Mathematics Program, under Grants DE-SC0025412 and DE-SC0026871.
Yanzhao Cao's research is partially supported by the U.S. Department of Energy, Office of Science,  under the grant number 	DE-SC0025649 and by the National Science Foundation under the grant number DMS-2446127. 
Zezhong Zhang acknowledges support from the U.S. National Science Foundation (NSF) under Grant DMS-2609189.

\bibliographystyle{plain}
\bibliography{references}

\end{document}